\documentclass[11pt,intlimits, draft,oneside]{amsart}

\usepackage[latin2]{inputenc}
\usepackage{latexsym}
\usepackage{amssymb}
\usepackage{amsthm}

\usepackage{amsmath}

\usepackage{enumerate}

\input{cyracc.def}
\newfont{\cyrfnt}{wncyi10 at 11pt}

\newtheorem{thm}{Theorem}[section]
\newtheorem{prop}[thm]{Proposition}
\newtheorem{lemma}[thm]{Lemma}

\theoremstyle{definition}
\newtheorem{defi}[thm]{Definition}

\newcommand{\R}{\mathbb{R}}

\newcommand{\C}{\mathbb{C}}
\newcommand{\N}{\mathbb{N}}
\newcommand{\U}{\mathcal{U}}

\renewcommand{\Re}{\mathrm{Re}\,}

\begin{document}

\title[Category theorems for stable operators]
{Category theorems for stable operators on Hilbert spaces}

\author{Tanja Eisner}
\address{Tanja Eisner \newline Mathematisches Institut, Universit\"{a}t T\"{u}bingen\newline Auf der Morgenstelle 10, D-72076, T\"{u}bingen, Germany}
\email{talo@fa.uni-tuebingen.de}


\renewcommand{\thefootnote}{}
\author{Andr\'{a}s Ser\'{e}ny} 
\footnote{The second author was supported by the DAAD-PPP-Hungary Grant, project number D/05/01422.} 
\address{Andr\'{a}s Ser\'{e}ny \newline
Department of Mathematics and its Applications, Central European
University,\newline N\'{a}dor utca  9, H-1051 Budapest, Hungary }
\email{sandris@elte.hu}


\begin{abstract}
We discuss the two closely related, but different concepts of  weak and almost weak stability for the powers of a contraction on a separable Hilbert space. 
Extending Halmos' and Rohlin's theorems in ergodic theory as a model, we show that the set of all weakly stable contractions is of first category while the set of all almost weakly stable contractions is of second category and is residual. 
Analogous statements for unitary and isometric operators are also proved.
\end{abstract}

\maketitle

\begin{center}
\emph{Dedicated to the inspiring work of the great mathematician Paul R.~Halmos.}
\end{center}

\vspace{0.1cm}


\section{Introduction}

In this paper we deal with the asymptotic behaviour of the powers of bounded linear operators on Hilbert spaces. 

The following classical decomposition theorem already yields some basic information. Here and later we denote the space of all bounded linear operators on a Hilbert space $H$ by $\mathcal{L}(H)$.
\begin{thm}\label{thm:JGdL}(Jacobs--Glicksberg--de Leeuw, see \cite[pp. 105--106]{krengel:1985})
Let $H$ be a Hilbert space and let $T\in\mathcal{L}(H)$ be power bounded, i.e., satisfying $\sup \Vert T^n\Vert <\infty$. Then $H=H_r\oplus H_s$, where
\begin{eqnarray*}
H_r&:=&\{x\in H:\ Tx=\gamma x \text{ for some } \gamma \text{ with } |\gamma|=1 \},\\
H_s&:=&\{x\in H:\ 0 \text{ is a weak accumulation point of } \{T^n x:\ n\in\N\}\}.
\end{eqnarray*}
\end{thm}
\noindent
(We note that the theorem is valid for every operator $T$ on an arbitrary Banach space $X$ such that each orbit $\{T^n x,\ n\in\N\}$ is weakly relatively compact in $X$. This holds for power bounded operators on reflexive Banach spaces.)   

The above theorem shows that in order to understand the asymptotic behaviour of an operator $T$ one needs to study its asymptotic behaviour on the subspace $H_s$. In particular, the question of convergence to $0$, called \emph{stability}, becomes central. 

In contrast to the well-established theory of uniform stability, \emph{weak stability}, i.e., the property 
\begin{equation*}
\lim_{n\to\infty} \langle T^n x,y \rangle =0\ \ \text{  for all } \ x,y\in H,
\end{equation*}
is much less understood even though it occurs naturally, e.g., see Sz.-Nagy, Foia\c{s} \cite[Chapter I]{sznagy/foias} and  Foguel \cite{foguel:1963}. In ergodic theory, see Halmos \cite[pp.~36--41]{halmos:1956} and Krengel \cite[p.~254]{krengel:1985}, it is called \emph{``mixing''} or \emph{``strong mixing''}. Indeed, Katok and Hasselblatt state on p. 748 of their monograph \cite{katok/hasselblatt:1995}:
\begin{quote}\it
``... It [strong mixing] is, however, one of those notions, that is easy and natural to define but very difficult to study...''
\end{quote}
The property satisfied by the restriction of $T$ to $H_s$ in Theorem \ref{thm:JGdL} is weaker, but much better understood. We will use the following terminology.
\begin{defi} 
A power bounded operator $T$ on a Hilbert space is called \emph{almost weakly stable} if $0$ is a weak accumulation point of every orbit $\{T^n x: n\in\N \}$.
\end{defi}
For power bounded operators on Hilbert spaces almost weak stability is, by Theorem \ref{thm:JGdL}, equivalent to the property ``no eigenvalues on the unit circle''. 

On bounded sets in separable Hilbert spaces the weak topology is metrisable and hence an operator $T$ is almost weakly stable if and only if for every $x\in H$ there exists a subsequence $\{n_k\}_{k=1}^\infty$ satisfying $\displaystyle\lim_{k\to\infty}\langle T^{n_k}x,y\rangle=0$ for every $y\in H$.  

%

Nagel \cite{nagel:1974}, Jones, Lin \cite{jones/lin:1976} and others have shown that almost weak stability for a power bounded operator $T$ implies convergence of $\langle T^n x,y\rangle$ to $0$ on a subsequence $\{n_k\}_{k=1}^\infty$ having asymptotic density $1$ for every $x,y\in H$. (The proof is based on the analogous statement from ergodic theory.) This explains the name ``almost weak stability''. However, on infinite-dimensional spaces this notion is not equivalent to weak stability, see Halmos \cite[pp.~77--80]{halmos:1956} for examples and e.g. \cite{EFNS} for a systematic discussion in the case of $C_0$-semigroups. 

Note that the concept of almost weak stability and its reformulations mentioned above are valid for power bounded operators on reflexive Banach spaces or, more generally, for operators having relatively compact orbits, as well.

It is our aim in this paper to show that (for operators on separable infinite-dimensional Hilbert spaces) the concepts of  weak stability and almost weak stability differ fundamentally. To be exact, we show that the sets of almost weakly and weakly stable operators have different Baire category in the complete metric spaces (with respect to an appropriate metric) of all unitary, isometric and contractive operators, see Theorems \ref{thm:unitary}, \ref{thm:isometry} and \ref{thm:contraction}, respectively. More precisely, the set of all weakly stable operators is of first category, while the set of all almost weakly stable operators is of second category and is even residual, i.e., its complement is of first category. In this sense, a typical operator of these classes is almost weakly but not weakly stable. These category results are analogous to classical theorems of Halmos and Rohlin for measure preserving transformations in ergodic theory, see Halmos \cite[pp.~77--80]{halmos:1956} or the original papers by Halmos \cite{halmos:1944} and Rohlin \cite{rohlin:1948}, where almost weak stability corresponds to the notion of weak mixing. Some constructions in our proof are similar to steps in Halmos' and Rohlin's proofs.

In the following we first treat the problem for unitary operators (Section \ref{section:unitary}), then for isometries (Section \ref{section:isometric}), and in Section \ref{section:contraction} for contractions.

\section{Unitary  case} \label{section:unitary}

Let $H$ be a separable infinite-dimensional Hilbert space. 
We denote the set of all unitary operators on $H$ by $\U$. 
The following density result for periodic operators is a first building block for our construction. 
(Here an operator $T$ is called periodic if there exists $n\in \N$ with $T^n=I$, the identity on $H$. 
The smallest such $n$ is called the period of the operator $T$.) 

\begin{prop}\label{prop:periodic} For every $n\in\N$ the set of all periodic unitary operators with period greater than $n$ is dense in $\U$ endowed with the operator norm topology.
\end{prop}

\begin{proof} Take $U\in \U$, $N\in \N$ and $\varepsilon>0$. By the spectral theorem $H$ is isomorphic to $L^2(\Omega, \mu)$ for some 
finite measure $\mu$ on a set $\Omega$ and $U$ is unitarily equivalent to a multiplication operator $\tilde{U}$ with 
\begin{equation*}
(\tilde{U}f)(\omega)=\varphi(\omega)f(\omega),\ \ \forall \omega\in\Omega,
\end{equation*} 
for some measurable $\varphi:\Omega \to \Gamma:=\{z\in\C:\ |z|=1\}$. 

We approximate the operator $\tilde{U}$ as follows. Consider the set
\begin{equation*}
\Gamma_N:=\{e^{2\pi i \cdot \frac{p}{q}}:\  p,q\in \N \text{ relatively prime },\ q> N\}
\end{equation*} 
which is dense in $\Gamma$. Take a finite set $\{\alpha_j\}_{j=1}^n \subset \Gamma_N$ such that $\arg(\alpha_{j-1})<\arg(\alpha_j)$ and $|\alpha_j - \alpha_{j-1}|<\varepsilon$ hold for all $2\leq j\leq n$. Define
\begin{equation*}
\psi(\omega):= \alpha_{j-1},\ \forall \omega\in \varphi^{-1}(\{z\in\Gamma:\ \arg(\alpha_{j-1})\leq \arg(z)< \arg(\alpha_j) \}). 
\end{equation*} 
Denote now by $\tilde{P}$ the multiplication operator with $\psi$. 
The operator $\tilde{P}$ is periodic with period greater than $N$. Moreover, 
\begin{equation*}
\Vert \tilde{U} - \tilde{P} \Vert = \sup_{\omega\in\Omega} |\varphi(\omega) - \psi(\omega)|\leq \varepsilon 
\end{equation*} 
holds and the proposition is proved.
\end{proof}
Before we present a second building block we need the following lemma.
\begin{lemma} \label{lemma:discrete-appr-of-I}
Let $H$ be a separable infinite-dimensional Hilbert space. Then there exists a sequence $\{T_n\}_{n=1}^\infty$ of almost weakly stable unitary operators satisfying $\displaystyle \lim_{n\to\infty}\Vert T_n - I\Vert=0$.  
\end{lemma}

\begin{proof}
By the isomorphism of all separable infinite-dimensional Hilbert spaces there exists a unitary operator $U:H\to L^2(\R)$, where $L^2(\R)$ is considered with the Lebesgue measure. 

Take $n\in\N$ and define $\tilde{T}_n$ on $L^2(\R)$ by
\begin{equation*}
(\tilde{T}_nf)(s):= e^{\frac{iq(s)}{n}}f(s), \ \ s\in\R,\ \ f\in L^2(\R),  
\end{equation*} 
where $q:\R\to [0,1]$ is strictly monotone.
Then all $\tilde{T}_n$ are almost weakly stable by the theorem of Jacobs--Glicksberg--de Leeuw and we have 
\begin{equation*}
\Vert \tilde{T}_n - I \Vert = \sup_{s\in\R} |e^{\frac{iq(s)}{n}} - 1| \leq |e^{\frac{i}{n}} - 1| \to 0 \ \text{  as  } \ n\to \infty.  
\end{equation*}
To finish the proof we only need to define $T_n:=U^*\tilde{T}_nU$ on $H$. 
\end{proof}

We now introduce the appropriate topology. It is called the \textit{strong* (operator) topology} and is induced by the family of seminorms $p_x(T):=\sqrt{\|Tx\|^2+\|T^*x\|^2}$, $x\in H$. We note that convergence in this topology corresponds to strong convergence of operators and their adjoints. For properties and further information on this topology we refer to Takesaki \cite[p.~68]{takesaki}. 

In the following we consider the space $\mathcal{U}$ of all unitary operators on $H$ endowed with the strong* operator topology. Note that $\U$ is a complete metric space with respect to the metric given by
\begin{equation*}
d(U,V):= \sum_{j=1}^\infty \frac{\|Ux_j -Vx_j\| + \|U^*x_j -V^*x_j\|}{2^j \|x_j\|}\quad  \text{for } U,V\in \U,
\end{equation*}
and $\{x_j\}_{j=1}^\infty$ some dense subset of $H\setminus \{0\}$. Further, by $\mathcal{S_U}$  we denote the set of all weakly stable unitary operators on $H$ and by $\mathcal{W_U}$ the set of all almost weakly stable unitary operators on $H$.

We now show the following density property for $\mathcal{W_U}$.

\begin{prop}\label{prop:almweakstab} The set  $\mathcal{W_U}$ of all almost weakly stable unitary operators is dense in $\U$.
\end{prop}

\begin{proof} By Proposition \ref{prop:periodic} it is enough to approximate periodic unitary operators by almost weakly stable unitary operators. Let $U$ be a periodic unitary operator and let $N$ be its period. Take $\varepsilon>0$, $n\in\N$ and $x_1,\ldots, x_n\in H\setminus\{0\}$. We have to find an almost weakly stable unitary operator $T$ with $\Vert Ux_j - Tx_j \Vert < \varepsilon$ and $\Vert U^*x_j - T^*x_j \Vert < \varepsilon$ for all $j=1,\ldots, n$. 

By $U^N=I$ and the spectral theorem, $\sigma(U)\subset \left\{1, e^{\frac{2\pi i}{N}},\ldots, e^{\frac{2\pi (N-1) i}{N}}\right\}$ and the orthogonal decomposition
\begin{equation}\label{Fix}
H=\ker(I-U) \oplus \ker(e^\frac{2\pi i}{N}I-U) \oplus \ldots \oplus \ker(e^\frac{2\pi (N-1) i}{N}I-U)
\end{equation} 
holds.

Assume first that $x_1,\ldots,x_n$ are orthogonal eigenvectors of $U$.

In order to use Lemma \ref{lemma:discrete-appr-of-I} we first construct a periodic unitary operator $S$ satisfying $Ux_j=Sx_j$ for all $j=1,\ldots, n$ and having  infinite-dimensional eigenspaces only. For this purpose define the n-dimensional $U$- and $U^*$-invariant subspace
$H_0 := \text{lin} \{x_j\}_{j=1}^n$
and the operator $S_0$ on $H_0$ as the restriction of $U$ to $H_0$. 
Decompose $H$ as an orthogonal sum
\begin{equation*}
H = \bigoplus_{k=0}^\infty H_k\quad \text{with } \text{dim}H_k=\text{dim}H_0\ \text{ for all } k\in\N . 
\end{equation*} 
Denote by $P_k$ an isomorphism from $H_k$ to $H_0$ for every $k$. Define now $S_k:=P_k^{-1}UP_k$ on each $H_k$ as a copy of $U|_{H_0}$ and consider $S:=\bigoplus_{k=0}^\infty S_k$ on $H$.

The operator $S$ is unitary and periodic with period being a divisor of $N$. So a decomposition analogous to (\ref{Fix}) is valid for $S$. Moreover, $Ux_j=Sx_j$ and $U^*x_j=S^*x_j$ hold for all $j=1,\ldots, n$ and the eigenspaces of $S$ are infinite dimensional. Denote by  $F_j$ the eigenspace of $S$ containing $x_j$ and by 
$\lambda_j$ the corresponding eigenvalue.  
By Lemma \ref{lemma:discrete-appr-of-I} for every $j=1,\ldots,n$ there exists an almost weakly stable unitary operator $T_j$ on 
$F_j$ 
satisfying $\Vert T_j - S_{|_{F_j}}  \Vert 
= \Vert T_j - \lambda_j I  \Vert 
< \varepsilon$. Finally, we define the desired operator $T$ as $T_j$ on $F_j$ for every $j=1\ldots,n$ and extend it linearly to $H$. 

Let now $x_1,\ldots,x_n\in H$ be arbitrary and take an orthonormal basis of eigenvalues $\{y_k\}_{k=1}^\infty$. Then there exists $K\in \N$ such that 
$x_j=\sum_{k=1}^K a_{jk} y_k + o_j$ with $\|o_j\|<\frac{\varepsilon}{4}$ for every $j=1,\ldots,n$. By the arguments above applied to $y_1,\ldots,y_K$ there is an almost weakly stable unitary operator $T$ with $\|Uy_k-Ty_k\|< \frac{\varepsilon}{4K M}$ and $\|U^*y_k-T^*y_k\|< \frac{\varepsilon}{4K M}$ for $M:=\max_{ k=1,\ldots,K, j=1,\ldots,n}{|a_{jk}|}$ and every $k=1,\ldots,K$. Therefore we obtain
$$\|Ux_j-Tx_j\|\leq \sum_{k=1}^K |a_{jk}|\|Uy_k-Ty_k\|+2\|o_j\|< \varepsilon$$ 
for every $j=1,\ldots,n$. Analogously, $\|U^*x_j-T^*x_j\|< \varepsilon$ holds for every $j=1,\ldots,n$, and the proposition is proved.
\end{proof}

We can now prove the following category theorem for weakly and almost weakly stable unitary operators.
To do so we extend the argument used in the proof of the corresponding category theorems for flows in ergodic theory (see Halmos \cite[pp.~77--80]{halmos:1956}).


\begin{thm}\label{thm:unitary}  
The set $\mathcal{S_U}$ of weakly stable unitary operators is of first category and the set  $\mathcal{W_U}$ of almost weakly stable unitary operators is residual in $\U$. 
\end{thm}

\begin{proof} First we prove that $\mathcal{S_U}$ is of first category in $\U$. Fix $x \in H$ with $\|x\|=1$ and consider 

\begin{equation*}
M_k:=\left\{U\in \U :\ |\langle U^k x,x\rangle| \leq \frac{1}{2} \right\}.
\end{equation*}
Note that all sets $M_k$ are closed. 

Let $U\in \U$ be weakly stable. Then there exists $n\in\N$ such that $U\in M_k$ for all $k\geq n$, i.e., $\displaystyle U\in\cap_{k\geq n} M_k$. So we obtain
\begin{equation}
\mathcal{S_U} \subset \bigcup_{n=1}^\infty N_n,
\end{equation}
where $N_n:=\cap_{k\geq n} M_k$. 
Since all $N_n$ are closed, it remains to show that $\U\setminus N_n$ is dense for every $n$. 

Fix $n\in\N$ and 
let $U$ be a periodic unitary operator. Then $U\notin M_k$ for some $k\geq n$ and therefore $U\notin N_n$. Since by Proposition \ref{prop:periodic} periodic unitary operators are dense in $\U$, $\mathcal{S_U}$ is of first category. 

To show that $\mathcal{W_U}$ is residual  we take a dense set $D=\{x_j\}_{j=1}^\infty$ of $H$ and define 
\begin{equation*}
W_{jkn}:=\left\{U\in \U :\ |\left<U^n x_j,x_j\right>| < \frac{1}{k} \right\}.
\end{equation*} 
All these sets are open. Therefore the sets 
$W_{jk}:=\cup_{n=1}^\infty W_{jkn}$ 
are also open.

We show that 
\begin{equation}\label{W}
\mathcal{W_U}=\bigcap_{j,k=1}^\infty W_{jk}
\end{equation}
holds. 

The inclusion ``$\subset$'' follows from the definition of almost weak stability. 
To prove the converse inclusion we take $U \in \mathcal{U}\setminus \mathcal{W_U}$.
Then there exists $x\in H$ with $\|x\|=1$ and $\varphi \in \R$ such that $Ux=e^{i\varphi} x$. 
Take now $x_j\in D$ with $\|x_j-x\|\leq \frac{1}{4}$. Then 
\begin{eqnarray*}
|\left<U^n x_j,x_j \right>| &=& |\left<U^n (x-x_j),x-x_j \right> + \left<U^n x,x \right> 
- \left<U^n x,x-x_j \right> - \left<U^n (x-x_j),x \right>| \\
&\geq& 1-\|x-x_j\|^2 -2\|x-x_j\|>\frac{1}{3}
\end{eqnarray*}  
holds for every $n\in\N$.
So 
$U\notin W_{j3}$ which implies 
$U\notin \cap_{j,k=1}^\infty W_{jk}$, and  
therefore 
(\ref{W}) holds. Moreover, all $W_{jk}$ are dense by Proposition \ref{prop:almweakstab}.
Hence $\mathcal{W_U}$ is residual as a countable intersection of open dense sets.
\end{proof}

\section{Isometric case} \label{section:isometric}

In this section we consider the space $\mathcal{I}$ of all isometries on $H$ endowed with the strong operator topology and prove analogous category results as in the previous section. We again assume $H$ to be separable and infinite-dimensional.
Note that $\mathcal{I}$ is a complete metric space with respect to the metric given by the formula 

\begin{equation*}
d(T,S):= \sum_{j=1}^\infty \frac{\|Tx_j -Sx_j\|}{2^j \|x_j\|}\quad \text{for } T,S \in \mathcal{I},
\end{equation*}
where $\{x_j\}_{j=1}^\infty$ is a fixed dense subset of $H\setminus \{0\}$. 

Further we denote by $\mathcal{S_I}$ the set of all weakly stable isometries on $H$ and by $\mathcal{W_I}$ the set of all almost weakly stable isometries on $H$.

The results in this section are based on the following classical theorem, called Wold decomposition of isometries on Hilbert spaces (see Sz.-Nagy, Foia{\c{s}} \cite[Theorem 1.1]{sznagy/foias}).
\begin{thm}\label{thm:Wold}
Let $V$ be an isometry on a Hilbert space $H$. Then $H$ can be decomposed into an orthogonal sum $H=H_0 \oplus H_1$ of $V$-invariant subspaces such that the restriction of $V$ on $H_0$ is unitary and the restriction of $V$ on $H_1$ is a unilateral shift, i.e., there exists a subspace $Y\subset H_1$ with $V^n Y\perp V^m Y$ for all $n\neq m$, $n,m\in\N$, such that $H_1=\oplus_{n=1}^\infty V^n Y$ holds.   
\end{thm}
 
We need the following easy lemma, see also Peller \cite{peller}.

\begin{lemma}\label{lemma:shift-appr}
Let $Y$ be a Hilbert space and let $R$ be the right shift on $H:=l^2(\N, Y)$. Then there exists a sequence $\{T_n\}_{n=1}^\infty$ of periodic unitary operators on $H$ converging strongly to $R$.
\end{lemma}

\begin{proof}
We define the operators $T_n$ by
\begin{equation*}
T_n(x_1, x_2, \ldots, x_n, \ldots):= (x_n, x_1, x_2, \ldots, x_{n-1}, x_{n+1},\ldots ).
\end{equation*}
Every $T_n$ is unitary and has period $n$. Moreover, for an arbitrary $x=(x_1,x_2,\ldots )\in H$ we have 
\begin{equation*}
\Vert T_n x - R x\Vert^2 = \|x_n\|^2 + \sum_{k=n}^\infty \|x_{k+1} - x_{k}\|^2 \to 0 \ \text{  as  } \ n\to\infty, 
\end{equation*}
and the lemma is proved.
\end{proof}

As a first application of the Wold decomposition we obtain the density of the periodic operators in  $\mathcal{I}$. (Note that periodic isometries are unitary.) 

\begin{prop}\label{prop:periodic-i} The set of all periodic isometries is dense in $\mathcal{I}$.
\end{prop}

\begin{proof} Let $V$ be an isometry on $H$. Then by Theorem \ref{thm:Wold} the orthogonal decomposition $H=H_0\oplus H_1$ holds, where the restricion $V_0$ on $H_0$ is unitary and the space $H_1$ is unitarily equivalent to $l^2(\N, Y)$. The restriction $V_1$ of $V$ on $H_1$ corresponds by this equivalence to the right shift operator on $l^2(\N, Y)$. By Proposition \ref{prop:periodic} and Lemma \ref{lemma:shift-appr} we can approximate both operators $V_0$ and $V_1$ by unitary periodic ones and the assertion follows.
\end{proof}

We  further obtain the density of the almost weakly stable operators in $\mathcal{I}$. 

\begin{prop}\label{prop:almweakstab-i} The set $\mathcal{W_I}$ of almost weakly stable isometries is dense in $\mathcal{I}$.
\end{prop}

\begin{proof} Let $V$ be an isometry on $H,\ H_0,\ H_1$ the orthogonal subspaces from Theorem \ref{thm:Wold} and $V_0$ and $V_1$ the corresponding restrictions of $V$. By Lemma \ref{lemma:shift-appr} the operator $V_1$ can be approximated by unitary operators on $H_1$. The assertion now follows from Proposition \ref{prop:almweakstab}.
\end{proof}

Using the same idea as in the proof of Theorem \ref{thm:unitary} one obtains with the help of Propositions \ref{prop:periodic-i} and \ref{prop:almweakstab-i} the following category result for weakly and almost weakly stable isometries. 

\begin{thm} \label{thm:isometry}
The set $\mathcal{S_I}$ of all weakly stable isometries is of first category and the set $\mathcal{W_I}$ of all almost weakly stable isometries is residual in $\mathcal{I}$. 
\end{thm}


\section{Contraction case} \label{section:contraction}

We now extend the category results in the previous sections to the case of contractive operators. The Hilbert space $H$ we take as before.

Let $\mathcal{C}$ denote the set of all contractions on $H$ endowed with the weak operator topology. Note that
 $\mathcal{C}$ is a complete metric space with respect to the metric given by the formula 

\begin{equation*}
d(T,S):= \sum_{i,j=1}^\infty \frac{|\left<Tx_i,x_j\right> -\left< Sx_i,x_j\right>|}{2^{i+j} \|x_i\| \|x_j\|}\quad \text{for }\ T,S \in \mathcal{C},
\end{equation*}
where $\{x_j\}_{j=1}^\infty$ is a fixed dense subset of $H\setminus \{0\}$.


By Takesaki \cite[p.~99]{takesaki}, the set of all unitary operators is dense in $\mathcal{C}$ (see also Peller \cite{peller} for a much stronger assertion). Combining this with Propositions \ref{prop:periodic} and \ref{prop:almweakstab} we have the following fact.

\begin{prop}\label{prop:density-contr}
The set of all periodic unitary operators and the set of all almost weakly stable unitary operators are both dense in  $\mathcal{C}$.
\end{prop} 


The 
property that weak convergence from below implies strong convergence is a key for the further results (cf. Halmos \cite[p. 14]{halmos:1967}). 

\begin{lemma}\label{w-imply-st}
Let $\{T_n\}_{n=1}^\infty$ be a sequence of linear operators on a Hilbert space $H$ converging weakly to a linear operator $S$. If $\|T_n x\|\leq \|Sx\|$ for every $x\in H$, then $\displaystyle \lim_{n\to\infty}T_n=S$ strongly.
\end{lemma}

\begin{proof}
For each $x\in H$ we have
\begin{eqnarray*}
\Vert T_n x - Sx \Vert^2 &=& \left< T_n x - Sx, T_n x - Sx \right> = \Vert Sx \Vert^2 + \Vert T_n x \Vert^2  - 2 \Re \left< T_n  x, S x \right> \\ &\leq& 2 \left< Sx,Sx \right>  - 2\Re\left< T_n x, Sx \right> = 2 \Re \left< (S - T_n) x, Sx \right> 
\to 0 \ \text{  as  } \ n\to\infty, 
\end{eqnarray*}
and the lemma is proved.
\end{proof}

We now state the category result for contractions. We note that its proof differs from the corresponding proofs in the previous sections. 

\begin{thm} \label{thm:contraction}
The set $\mathcal{S_C}$ of all weakly stable contractions is of first category and the set $\mathcal{W_C}$ of all almost weakly stable contractions is residual in $\mathcal{C}$. 
\end{thm}

\begin{proof}
To prove the first statement we fix $x\in H$, $\Vert x \Vert =1$, and define as before the sets

\begin{equation*}
N_n:=\left\{T\in \mathcal{C} :\ |\langle T^k x,x\rangle | \leq \frac{1}{2} \text{  for all } k\geq n \ \right\}.
\end{equation*}
%
Let $T\in \mathcal{C}$ be weakly stable. Then there exists $n\in\N$ such that $\displaystyle T\in N_n$, and we obtain
\begin{equation}
\mathcal{S_C} \subset \bigcup_{n=1}^\infty N_n.
\end{equation}

\noindent 
It remains to show that the sets $N_n$ are nowhere dense. 
%
Fix $n\in\N$ and 
let $U$ be a periodic unitary operator. We show that $U$ does not belong to the closure of $N_n$. Assume the opposite, i.e., that  there exists a sequence $\{T_k\}_{k=1}^\infty \subset N_n$ satisfying $\lim_{k\to\infty}T_k=U$ weakly. Then by Lemma \ref{w-imply-st} $\lim_{k\to\infty}T_k=U$ strongly and therefore $U\in N_n$ by the definition of $N_n$. This contradicts the periodicity of $U$. By the density of the set of unitary periodic operators in $\mathcal{C}$ we obtain that $N_n$ is nowhere dense and therefore $\mathcal{S_C}$ is of first category.


To show the residuality of $\mathcal{W_C}$  we again take a dense subset $D=\{x_j\}_{j=1}^\infty$ of $H$ and define 
\begin{equation*}
W_{jk}:=\left\{T\in \mathcal{C} :\ |\langle T^n x_j,x_j\rangle | < \frac{1}{k} \text{ for some } n\in\N \right\}.
\end{equation*} 

\noindent  As in the proof of Theorem \ref{thm:unitary} the equality 
\begin{equation}\label{W_C}
\mathcal{W_C}=\bigcap_{j,k=1}^\infty W_{jk}
\end{equation}
holds. 

Fix $j,k\in\N$. We have to show that the complement $W^c_{jk}$  of $W_{jk}$ is nowhere dense. We note that 
$$W^c_{jk}=\left\{T\in \mathcal{C} :\ |\left<T^n x_j,x_j\right>| \geq \frac{1}{k} \text{ for all } n\in\N \right\}.$$ 
Let $U$ be a unitary almost weakly stable operator. Assume that there exists a sequence $\{T_m\}_{m=1}^\infty \subset W^c_{jk}$ satisfying $\lim_{m\to\infty}T_m= U$ weakly. Then, by Lemma \ref{w-imply-st}, $\lim_{m\to\infty}T_m= U$ strongly and therefore $U\in W^c_{jk}$. This contradicts the almost weak stability of $U$. Therefore the set of all unitary almost weakly stable operators does not intersect the closure of $W^c_{jk}$. By Proposition \ref{prop:density-contr} all sets $W^c_{jk}$ are nowhere dense and therefore $\mathcal{W_C}$ is residual. 
\end{proof}


\noindent {\bf Final remark.}
It is not clear in which reflexive Banach spaces results similar to Theorems \ref{thm:isometry} and \ref{thm:contraction} hold. Further, it is not clear whether  in Theorem \ref{thm:contraction} one can replace contractions by power bounded operators.

\vspace{0.3cm}

\noindent {\bf Acknowledgement.} 
The authors are very grateful to Andr\'as B\'atkai and Rainer Nagel for valuable comments and interesting discussions. We also express our sincere thanks to the editor L\'aszl\'o K\'erchy and the referee, whose suggestions have considerably improved the paper. 

\parindent0pt

\end{document}